\input amstex
\documentstyle{amsppt}
\document

\magnification 1100

\def\gen{{\frak{g}}}

\def\len{\frak{l}}

\def\sen{\frak{s}}

\def\a{{\alpha}}
\def\o{{\omega}}

\def\l{{\lambda}}
\def\g{{\gamma}}

\def\b{{\beta}}
\def\eps{{\varepsilon}}

\def\1b{{\bold 1}}

\def\cb{{\bold c}}
\def\db{{\bold d}}
\def\eb{{\bold e}}
\def\fb{{\bold f}}
\def\gb{{\bold g}}

\def\ib{{\bold i}}
\def\jb{{\bold j}}

\def\Bb{{\bold B}}

\def\Mb{{\bold M}}
\def\Nb{{\bold N}}
\def\Hb{{\bold H}}

\def\Rb{{\bold R}}
\def\Sb{{\bold S}}
\def\Tb{{\bold T}}

\def\Ub{{\bold U}}
\def\Vb{{\bold V}}
\def\Wb{{\bold W}}

\def\X{{\roman X}}

\def\W{{\roman W}}

\def\T{{\roman T}}

\def\wt{{\text{wt}}}

\def\sm{{\text{sm}}}

\def\Id{\text{Id}\,}

\def\simto{\,{\buildrel\sim\over\to}\,}

\def\wt{\roman{wt}}

\def\Im{\text{Im}\,}

\def\wt{\text{wt}}

\def\GL{{\text{GL}}}

\def\T{{\text{T}}}

\def\Otimes{\ts\bigotimes}
\def\Sum{\ts\sum}

\def\AA{{\Bbb A}}

\def\KK{{\Bbb K}}

\def\NN{{\Bbb N}}

\def\QQ{{\Bbb Q}}
\def\RR{{\Bbb R}}

\def\ZZ{{\Bbb Z}}

\def\Ac{{\Cal A}}
\def\Bc{{\Cal B}}

\def\Fc{{\Cal F}}

\def\wt{{\text{wt}}}
\def\wta{{\text{wt}^a}}
\def\and{{\ \text{and}\ }}

\def\ts{\textstyle}

\def\qed{\hfill $\sqcap \hskip-6.5pt \sqcup$}        
\overfullrule=0pt                                    

\def\gld{{{\gen\len_d}}}

\def\glp{{{\gen\len_p}}}
\def\gln{{{\gen\len_n}}}
\def\slp{{{\sen\len_p}}}
\def\lra{{{\longrightarrow}}}

\def\u1{{\underline 1}}

\def\la{{\langle}}
\def\ra{{\rangle}}

\newdimen\Squaresize\Squaresize=14pt
\newdimen\Thickness\Thickness=0.5pt
\def\Square#1{\hbox{\vrule width\Thickness
	      \alphaox to \Squaresize{\hrule height \Thickness\vss
	      \hbox to \Squaresize{\hss#1\hss}
	      \vss\hrule height\Thickness}
	      \unskip\vrule width \Thickness}
	      \kern-\Thickness}
\def\Vsquare#1{\alphaox{\Square{$#1$}}\kern-\Thickness}

\document

\topmatter
\title Periodic modules and quantum groups\endtitle
\rightheadtext{Periodic modules and quantum groups}
\leftheadtext{Michela Varagnolo}
\author Michela Varagnolo \endauthor
\address D\'epartement de Math\'ematique, Universit\'e de Cergy-Pontoise,
2 av. A. Chauvin, BP 222, 95302 Cergy-Pontoise cedex, France\endaddress
\email michela.varagnolo\@math.u-cergy.fr\endemail

\thanks
\endthanks

\abstract  
We prove that the elements $A_\leq$ defined by Lusztig
in a completion of the periodic
module actually live in the periodic module (in the type $A$ case).
In order to prove this, we compare, using the Schur duality,
these elements with Kashiwara canonical basis of an integrable module.
\endabstract
\endtopmatter
\document
\head 0. Introduction\endhead
Fix two positive integers $d$ and $p.$
Let $V'\supset R\supset R^+\supset I=\{\a_1,...\a_{d-1}\}$
be a $\RR$-vector space of dimension $d-1$, 
a root system of type $A_{d-1},$ a system of positive
roots and the corresponding set of simple roots. Fix a partition $\cb=(c_1,
...c_\ell)$ of $d$ and set $I_\cb=\{\a_i\,|\,i\not= c_1+\cdots c_a,\forall a\}.$
Let $\Fc$ be the set of all reflection hyperplanes $F_{\a^\vee,n}=
\{v\in V'\,|\, (\a^\vee:v)=np\}$, where $\a\in R^+$ and $n\in\ZZ.$
Set $\Fc_\cb=\{F_{\a^\vee,n}\,|\, \a\in R^+\cap I_\cb\}\subset \Fc.$
We call alcoves the connected components of $V'\setminus\cup_{F\in \Fc}F.$
Let $A'_+$ be the only alcove contained in the dominant Weyl chamber having
0 in its closure. 
Let $\Ac_\cb$ be the set of  alcoves contained in the connected 
component of $V'\setminus\cup_{F\in \Fc_\cb}F$ which contains $A'_+.$
Finally, let $\Mb_\cb$ be the $\ZZ[q,q^{-1}]$-span of  the alcoves 
in $\Ac_\cb.$ It is a $\Hb'\times \Rb'_\cb$-module, where $\Hb'$
is the affine Hecke algebra associated to $R$, and $\Rb'_\cb$ is
a Laurent polynomial algebra over $\ell$ variables.
In \cite{L5, 9.17(a)} Lusztig defines,
for any $A\in\Ac_\cb,$ some element $A_\leq$ in a completion of $\Mb_\cb$
and he conjectures in \cite{L5, 12.7} that $A_\leq\in \Mb_\cb$ 
(the conjecture is stated for all types).
If $\cb=(1^d),$ the conjecture follows from results in \cite{L1}.
In this paper we prove it for all $\cb$ in Theorem 5.5.

In \cite{L5, 17.3} two conjectural multiplicity formulas for modules
of a simple Lie algebra over an algebraically closed
field of characteristic $p$ are given. They are formulated via some 
 element $A_\geq$ and it involves  polynomials $\pi_{AB},$
 the coefficients of the expansion of $A_\geq$ in the basis of the alcoves,
\cite{L5, 9.17(c),(d)}. 
Using Theorem 5.5 and \cite{L5, Proposition 12.9}
we get that for all $A,$ there exist only finitely many non-zero
polynomials $\pi_{AB}.$

The main tools of the proof are the Schur duality (in its quantum affine
version) and a theorem of Kashiwara yielding a canonical basis
for some integrable module $\Nb$ of the quantum affine loop 
algebra of $\glp$ (here $p>d$). 
More precisely we write down a $\Hb'\times \Rb'_\cb$-isomorphism
which sends a $\Hb'\times \Rb'_\cb$-cyclic submodule of 
$\Mb_\cb$ into a weight space of $\Nb$ and takes $\{A_\leq\}$
into Kashiwara's canonical basis.

This work is motivated by \cite{VV3}. In {\it loc.cit.} the authors give a
geometric construction, via quiver varieties, of
the (signed) canonical basis of the module $\Nb(\l_\cb).$ This construction
is very similar to Lusztig's conjectural construction of the (signed) canonical
basis in the K-theory of the Springer fiber
(see \cite{L4, L5}, where it is also explained the link between this K-theory
object and the periodic module).
It is known that, for type $A$, the Springer fibers are quiver varieties.
Probably it is possible to check that both signed bases coincide via their
geometric characterization.
In this paper we give an elementary algebraic proof.

{\it Aknowledgements.} The idea to compare the elements $A_\leq$ with the 
Kashiwara's canonical basis is due to  E. Vasserot.
I thank him to share it with me.
This work started during my stay at the M.S.R.I. in March 2002. I thank
the organizers of the program "Infinite-dimensional algebras and Mathematical
Physics" for inviting me.

\head 1. Notations\endhead

\subhead 1.1\endsubhead
Let $\T\subset\GL_d$ be a maximal torus.
Let $\X=\X(\T)$ be the caracter group of $\T,$
and let $\Rb(\T)$ be its representation ring.
The dual group of one parameter subgroups of
$\T$ is denoted by $\X^\vee$.
Set $(\,:\,)$ equal to the canonical pairing
$\X\times\X^\vee\to\ZZ$.
Set $R\subset\X$ equal to the set of root of $\gld$
with respect to $\T$.
Choose a positive system $R^+\subset R$.
Let $I\subset R^+$ be the set of simple roots and let $Q\subset\X$
be the root lattice.
We fix basis elements $\eps_1,...,\eps_d\in\X$ such that
$I=\{\a_i=\eps_i-\eps_{i+1}\,|\,i=1,2,...,d-1\}.$
For $i=1,...,d$ set $\o_i=\eps_1+\cdots+\eps_i\in\X$ with $i=1,...,d$.
Let $\{\o^\vee_i\}, \{\eps^\vee_i\}\subset\X^\vee$ be the dual bases.
Put $\theta^\vee=\a_1^\vee+\cdots+\a^\vee_{d-1},$ and
$\a^\vee_i=\eps^\vee_i-\eps^\vee_{i+1}$ with $i=1,2,...,d-1.$
Let $\leq$ be the partial order on $\X$
such that $\mu\leq\nu$ if and only if
$\nu-\mu\in\Sum_{i=1}^{d-1}\NN\,\a_i$.
\subhead 1.2\endsubhead
Let $\W^f$ be the Weyl group of $\gld$.
Let $\W=\W^f\ltimes\X$,
$\W'=\W^f\ltimes Q,$
be the extended affine Weyl group and the affine Weyl group.
For any $w\in\W^f$, $\l\in\X$ we use the following notations :
$w=(w,0)$, $\tau_\l=(0,\l).$

A composition of $d$ is a $\ell$-uple of non-negative integers
$\fb=(f_1,f_2,...f_\ell)$ such that $\sum_i f_i=d$, for some $\ell.$
It is a partition if and only if $f_1\geq f_2\geq\cdots\geq f_\ell\geq 0$.
For any composition $\fb$ of $d$, set 
$I_\fb=\{\a_i\in I\,|\, i\neq f_1+f_2+\cdots +f_a,\ \forall a\}.$
Let $\W_\fb$
be the corresponding parabolic subgroup of $\W^f,$
i.e. $\W_\fb=\langle s_i\,|\,\a_i\in I_\fb\rangle.$ 
Put $\W'_\fb=\W_\fb\ltimes \ZZ I_\fb.$ 
Let $\W^\fb$(resp. ${}^\fb\W$)$\subset\W^f$ be the set of all elements
$w$ having a minimal length in the coset $w\W_\fb$ (resp. $\W_\fb w$).
Let $w_d$, $w_\fb$ be the longest elements of
$\W^f,$ $\W_\fb.$ 
Set $\nu_d=\ell(w_d), \nu_\fb=\ell(w_\fb)$.

Let $s_1,s_2,...,s_d$ be the simple affine reflections in $\W'$
such that $s_i=s_{\a_i^\vee}$ if $i\neq d$ and 
$s_d=\tau_\theta\, s_{\theta^\vee},$
where for any root $\a\in R$, $s_{\a^\vee}$
is the corresponding reflection.
Let $\pi$ denote the element $\tau_{\o_1}s_1s_2\cdots s_{d-1}.$
The cyclic group of infinite order $\la\pi\ra$ is isomorphic to
$\X/ Q$ and $\W$ is isomorphic to a semi-direct product
$\la\pi\ra\ltimes\W'.$
\subhead 1.3\endsubhead
Fix a positive number, say $p$. Set $V=\X\otimes_\ZZ\RR,$
$$A_+=\{\g\in V\,|\,p>(\g:\a^\vee)>0,\,\forall\a\in R^+\},$$
$$\X_p=\{\g\in\X\,|\,p\geq(\g:\eps_i^\vee)> 0,\,\forall i=1,...d\}.$$
The group $\W$ acts on $\X$ as follows
$$\matrix
\g\cdot s_i=\g-(\g:\a_i^\vee)\,\a_i\ \hfill&
\ \text{if}\ i=1,2,...,d-1,\hfill\cr
\g\cdot\tau_\nu=\g-p\nu\ \hfill&\ \text{if}\ \nu\in\X.\hfill
\endmatrix$$
Note that

- the set $\bar A_+\cap \X_p$ is a fundamental domain for the
right action of $\W$ on $\X$ ($\bar A_+$ is the closure of $A_+$ in $V$),

- for all $\mu\in \bar A_+\cap \X_p$,  $\W_\eb$ is the isotropy group of $\mu,$
if $\eb$ is associated to $\mu$ as in (1.4.1).
\subhead 1.4\endsubhead
{\bf Conventions.}

(1) In the future we will use both weights of $\glp$ and $\gld$.
The fundamental weights and the simple roots of $\glp$ are denoted
by $\Lambda_i$, $i=1,...,p$ and $\b_i=\epsilon_i-\epsilon_{i+1}$,
$i=1,...,p-1,$ respectively.
The affine simple roots of $\glp$ are $\b_1,...,\b_p.$ Set 
$$\delta=\b_1+\cdots +\b_p,\qquad\tilde\X=\Sum_{a=1}^p\ZZ
\epsilon_a,\qquad\tilde\X^a=\ZZ\delta\oplus\tilde\X.$$

(2) From now on, 
$\tilde\l=\Sum_{a=1}^{p}d_a\epsilon_a$ will be a fixed  (dominant) 
weight with  $d_1\geq d_2\geq\cdots\geq  d_p\geq 0,$ $\sum_{a=1}^pd_i=d.$
Set $\db= (d_p,...,d_1)$ and let 
$\cb=(c_1,c_2,...,c_\ell)$  be the partition dual to $\db$, that
is $c_i=|\{j\,|\,d_j\geq i\}|$. Then
$\cb=(p^{\ell_p}...2^{\ell_2}1^{\ell_1}),$  
for some positive integers $\ell_1,...,\ell_p$ such that
$\ell=\Sum_{a=1}^p\ell_a.$ We have
$$\tilde\l=\Sum_{a=1}^{p}\ell_a\Lambda_a=
\Sum_{a=1}^\ell \Lambda_{c_a}.$$

(3) Given  $\mu=\sum_{i=1}^{d}b_i\eps_i\in \X_p,$ we set
$$\tilde\mu=\Sum_{a=1}^pe_a\epsilon_a\in\tilde\X
\ \text{ with  }\ e_a=|\{i\,|\,b_i=a\}|,\qquad
\eb=(e_p,...,e_2,e_1).
\leqno(1.4.1)$$
Set 
$$\Omega=\{\Sum_{a=1}^p e_a\epsilon_a\in\tilde\X\,|\, 
e_a\geq 0, d=\Sum_ae_a\},
\Omega_\sm=\{ \Sum_{a=1}^p e_a\epsilon_a\in\Omega\,|\,
\eb \text{ is small}\},$$
where $\eb$ is  as in (1.4.1). 
The composition $\eb$ is said  small if $e_a=0,1$ for all $a.$
The map $\mu\mapsto \tilde\mu$ restricts to a bijection $\bar{A}_+\cap\X_p\to
\Omega$ which takes $A_+\cap\X_p$ onto $\Omega_\sm.$
From now on if $\tilde\mu\in \Omega$ (resp. $\Omega_\sm$), then $\mu$
denotes the
unique preimage  in $\bar A_+\cap\X_p$ (resp. $A_+\cap\X_p$). 

\head 2. $q$-Notations\endhead

\subhead 2.1\endsubhead
Put 
$\AA=\QQ[q,q^{-1}],\KK=\QQ(q), 
\AA_0=\{f\in\KK\,|f\text{ is without pole at }q=0\}.$
Set 
$$\Rb=\ZZ[\X]\otimes\AA,\quad\Rb_\cb=\ZZ[\X/\ZZ I_\cb]
\otimes \AA,\quad
\Rb'=\ZZ[Q]\otimes \AA,\quad 
\Rb'_\cb=\ZZ[Q/\ZZ I_\cb]\otimes \AA.$$
We identify $\Rb$ with the polynomial ring $\AA[x_1^{\pm 1},
x_2^{\pm 1},...,x_d^{\pm 1}],$ via the $\AA$-algebra isomorphism
which takes  $\eps_i\otimes 1$ to $x_i$.
For all $\g=\sum_i m_i\eps_i\in\X,$ put $x_\g=\prod_ix_i^{m_i}.$
Then $\Rb_\cb$ is identified with the quotient of $\Rb$ by the
ideal generated by the relations $x_{\a_i}=1,$ for all $\a_i\in I_\cb,$ and
$\Rb'$ with the sub-$\AA$-algebra of $\Rb$
generated by $\{x_{\a_i}^{\pm 1}\,|\,\a_i\in I\}$. We denote by 
$x_{\g+\ZZ I_\cb}$ the image of $x_\g$ in $\Rb_\cb.$
Set
$\a(c_i)=\sum_{k=1}^{c_i}(c_i+1-2k)\eps_{c_1+\cdots+c_{i-1}+k},$
and $\a_\cb=\sum_{i=1}^\ell \a(c_i)\in Q.$
Define an $\AA$-algebra homomorphism $\psi : \Rb\to \Rb_\cb$
by the formula $\psi(x_\g)=q^{(\g:\a_\cb)}x_{\g+\ZZ I_\cb}.$
The map  $\psi$ restricts to 
$\psi: \Rb'\to\Rb'_\cb,$ cf. \cite{L5, 8.3}.
We have
$$\psi(x_{\a_i})=q^2\qquad\text{for all }\a_i\in I_\cb.
\leqno(2.1.1)$$
\subhead 2.2\endsubhead
The quantum loop algebra of $\glp$ is the
$\KK$-algebra $\Ub_{\KK}$ generated by
elements $e_a,f_a,l_a^{\pm 1}$,
$a=1,2,...,p$,
modulo the following defining relations
$$l_al_a^{-1}=1=l^{-1}_al_a,\quad l_al_b=l_bl_a,$$
$$l_be_al_b^{-1}=q^{\delta(a= b)-\delta(a+1\equiv b)}e_a,\quad
l_bf_al_b^{-1}=q^{\delta(a+1\equiv b)-\delta(a= b)}f_a,$$
$$[e_a,f_b]=
\delta(a,b){{k_a-k_a^{-1}}\over q-q^{-1}},$$
$$\sum_{p=0}^m(-1)^{^p}
\left[\matrix m\cr p\endmatrix\right]
e_a^pe_be_a^{m-p}=
\sum_{p=0}^m(-1)^{^p}
\left[\matrix m\cr p\endmatrix\right]
f_a^pf_bf_a^{m-p}
=0.$$
In the last identity  $a\neq b$, $m=2$  if $a-b=\pm 1$ and 1 else.
We have set $[n]=q^{1-n}+q^{3-n}+...+q^{n-1}$ if $n\geq 0$,
$[n]!=[n][n-1]...[2]$, and
$$\left[\matrix m\cr p\endmatrix\right]=
{[m]!\over[p]![m-p]!}.$$  
We have also set $k_a=l_al_{a+1}^{-1},$ for all $a=1,2,...,p-1$ and
$k_p=l_pl_1^{-1}.$ Finally $a\equiv b$ means $a-b\in p\ZZ.$
Let $\Delta$ be the coproduct of $\Ub_{\KK}$, defined as follows
$$\Delta(e_a)=1\otimes e_a+e_a\otimes k_a^{-1},\quad
\Delta(f_a)= k_a\otimes f_a+f_a\otimes 1,\quad
\Delta(l_a)=l_a\otimes l_a.$$

The algebra $\Ub_\KK$  has also a presentation in terms of 
Drinfeld generators $x_{ir}^\pm,h_{is}^\pm,l_a^{\pm 1},$
where $r\in\ZZ,s\in\ZZ\setminus\{0\}, i=1,2,...,p-1,a=1,2,...,p.$
We normalize these generators as in \cite{B}
(with $T_i=T_{i,1}^{''}$ in the  notations of \cite{L2}). The presentation
depends on a choice of a 
function $o:\{1,...p-1\}\to\{\pm 1\}$ such that $o(i\pm 1)=-o(i).$

Let $\dot\Ub_{\KK}$ be the modified algebra of $\Ub_{\KK}$. It is a 
$\KK$-algebra without unity generated by $\pi_\ib e_a, \pi_\ib f_a,\pi_\ib$,
$\ib\in\NN^p, a=1,2,...p,$ with $\pi_\ib\pi_\jb=\delta_{\ib\jb}\pi_\ib$,
see \cite{L2, 23.1} for a precise definition.
Let $\Ub$, $\dot\Ub$ be the $\AA$-forms of
$\Ub_{\KK}$, $\dot\Ub_{\KK}.$ Then  $\Ub$ is the $\AA$-subalgebra
of $\Ub_\KK$ generated by $e_a^{(n)},f_a^{(n)},l_a^{\pm 1},$ and
$\dot\Ub$ is the $\AA$-subalgebra of $\dot\Ub_\KK$
generated by $\pi_\ib e_a^{(n)},\pi_\ib f_a^{(n)}$. 

Let $x\mapsto\overline{x}$ be the ring homomorphism of $\Ub$ such that
$\bar q=q^{-1}$, $\bar e_a=e_a$, $\bar f_a= f_a$, $\bar l_a=l_a^{-1}$.
We will still denote by $x\mapsto\overline{x}$ the involution on 
$\dot{\Ub}$ which is also characterized by $\overline{\pi_\ib}=\pi_\ib.$
\vskip2mm
\noindent{\bf Remark.}
Kashiwara's  canonical bases are defined for the quantum loop 
algebra of $\slp$. In this paper 
we will use the $\glp$-version. 
It is irrelevant since
an integrable  $\Ub$-module is a module for 
the quantum loop algebra of $\slp$ together with
an extra grading coming from the action of the $l_a$'s.
\subhead 2.3\endsubhead
Let $\Hb$ be the affine Hecke algebra of $\GL_d$.
Recall that $\Hb$ is the unital associative
$\AA$-algebra generated by elements
$t_1$,...$t_{d-1}$,
$x_1^{\pm 1}$,...$x_d^{\pm 1}$, modulo the relations
$$\matrix
(t_i+q^{-1})(t_i-q)=0,\hfill
&\hfill\cr
t_it_{i+1}t_i=t_{i+1}t_it_{i+1},\hfill
&|i-j|\geq 2\Rightarrow t_it_j=t_jt_i,\hfill\cr
x_ix_i^{-1}=1=x_i^{-1}x_i,\hfill
&x_ix_j=x_jx_i,\hfill\cr
t_ix_it_i=x_{i+1},\hfill
&j\neq i,i+1\Rightarrow t_ix_j=x_jt_i.\hfill
\endmatrix$$
For any $w\in\W$ let $t_w$ be the corresponding element in $\Hb.$
To simplify we set $\pi=t_\pi\,(=
x_1 t_1 t_2\cdots t_{d-1}).$ Set also  $t_d=\pi t_{d-1}\pi^{-1}$.
Let $\Hb',\Hb^f\subset\Hb$ be the $\AA$-subalgebras generated
by $t_1,...t_d$ and $t_1,...t_{d-1}$ respectively.
We will identify, in the obvious way,
$\Rb$ with the sub-$\AA$-algebra of
$\Hb$ generated by the $x^{\pm 1}_i$'s.
Note that $\Hb'$ is the sub-$\AA$-algebra of $\Hb$
generated by $t_1,...t_{d-1}$ and by $\Rb'.$
For future references, we point out 
that $x_\a$ is equal to the element $\theta_{-\a}$ 
in Lusztig's papers.

Let $t\mapsto\bar t$ be the ring homomorphism of $\Hb$ such that
$\bar q=q^{-1}$, $\bar t_i=t_i^{-1}$, $\bar \pi=\pi$, $i=1,...,d.$
\subhead 2.4\endsubhead
For any  composition $\fb$ of $d$, 
let $\Hb_\fb\subseteq\Hb^f$
be the subalgebra generated by the elements $t_i$
with $\a_i\in I_\fb.$
Similarly, let $\Hb'_\fb\subseteq\Hb'$ be the subalgebra
generated by $\Hb_\fb$ and the elements $x^{\pm 1}_{\a_i}$
with $i\in I_\fb$.
Set also 
$$\rho_\fb=\sum_{w\in\W_\fb}q^{\ell(w)} t_w,
\qquad m_\fb=\sum_{w\in\W_\fb} q^{2\ell(w)}.$$
We get 
$\rho_\fb^2= m_\fb\rho_\fb$ and $\overline{\rho_\fb}=q^{-2\nu_\fb}\rho_\fb.$
\subhead 2.5 \endsubhead
Set $\AA^-$ equal to the unique left representation of $\Hb'$
on $\AA$ taking the elements $t_i$ to $-q^{-1}\Id$, for all $i=1,...d.$
Note that, in particular,
$$x_{\a_i}\cdot 1= q^2\cdot 1\qquad \text{for all }\a_i\in I
\leqno(2.5.1).$$
\subhead 2.6\endsubhead
For any pair of compositions $\fb,\fb'$ of $d$ with at most $p$ parts,
let  $\Hb_{\fb,\fb'}$ be the 
$\AA$-linear span of the elements $t_m=\sum_{w\in m}q^{\ell(w)}t_w$,
where $m\in \W_\fb\backslash\W/\W_{\fb'}.$
Let $\Sb_d=\bigoplus_{\fb,\fb'}\Hb_{\fb,\fb'}$ the $q$-affine Schur algebra.
The product in $\Sb_d$ is defined 
by setting 
$t_m\star t_n=\delta(\fb',\gb)m_{\fb'}^{-1}t_m t_n$ for all
$m\in\W_\fb\backslash\W/\W_{\fb'},$
$n\in\W_\gb\backslash\W/\W_{\gb'}.$
Consider the  algebra homomorphism $\chi_d:\dot{\Ub}\to\Sb_d,$
defined in \cite{SV, Proposition 2.4}. In particular it takes 
$\pi_\ib$ to $\rho_{\ib^\sharp}$, 
if $\ib^\sharp$ is a composition of $d$, and
$\pi_\ib\mapsto 0$ otherwise.
Here  $\ib^\sharp=(i_1,i_2,...)$ and
$i_r=|\ib^{-1}(r)|.$  Let
$\chi_{\fb,\fb'}:\dot{\Ub}\to\Hb_{\fb,\fb'}$
be the composition of $\chi_d$ with the
canonical projection $\Sb_d\to\Hb_{\fb,\fb'}.$
The $\AA$-module $\Sb_d$ has a basis, $\Bc=\{b_\sen\,|\,\sen\in\Ac_d\}$,
where $\Ac_d$ is the following subset of the set of $\ZZ\times\ZZ$ matrices
with entries in $\NN,$
$$\Ac_d=\{\sen=(s_{ij})_{i,j\in\ZZ}\,|\, s_{i+p,j+p}=s_{ij},
\,\sum_{i\in\ZZ}\sum_{j=1}^ps_{ij}=d\},$$
see \cite{SV, 3.2}.
Then \cite{L6, Theorem 8.2} states that 
$\Bc^{ap}=\{b_\sen\,|\, \sen\in\Ac^{ap}_d\}$ is a $\AA$-basis of $\Im(\chi_d),$
with
$$\Ac_d^{ap}=\{\sen\in\Ac_d\,|\, \forall j\in\ZZ\setminus\{0\},\exists
i\in\ZZ\and s_{i,i+j}=0\}.$$
For a future use let us mention the following.
\proclaim{Claim} 
If $p>d$ and $\fb$ is small, the map $\chi_{\fb,\fb'}$ is 
surjective for all $\fb'.$ If $p=d$ and $\fb,\fb'$ are small, then 
$\Bc\cap\Hb_{\fb,\fb'}=(\Bc^{ap}\cap\Hb_{\fb,\fb'})\cup
\{\rho_\fb\pi^r\,|\, r\in\ZZ\}.$
\endproclaim
\demo{Proof} Set $\fb=(f_1,...,f_p),\fb'=(f'_1,...,f'_p).$ Then 
$$b_\sen\in\Hb_{\fb,\fb'}\Leftrightarrow\sen\in\Ac_{\fb,\fb'}=
\{\sen\in\Ac_d\,|\, \sum_{j\in\ZZ}s_{ij}=f_i,\sum_{i\in\ZZ}s_{ij}=f'_j,
\forall i,j\in\{1,2,...p\}\}.$$
Suppose that $p\geq d.$ If $\fb$ is small,  a direct computation gives
$$\Ac_{\fb,\fb'}\cap\Ac_d^{ap}\not=\Ac_{\fb,\fb'}\Leftrightarrow
\fb'\text{ is small}\,\and\, p=d.$$
This proves the first part of the claim. 
If  $p=d$ and $\fb,\fb'$ are both small, we have
$$\Ac_{\fb,\fb'}\setminus(\Ac_{\fb,\fb'}\cap\Ac_d^{ap})=
\{\sen(r)\,|\,r\in\ZZ\setminus\{0\}\},$$
with $\sen(r)=(s_{ij})_{i,j\in\ZZ}$ such that $
s_{ij}\not=0\Leftrightarrow j-i=r.$
Moreover $b_{\sen(r)}=\rho_\fb\pi^r,$ see \cite{SV, 3.2}.
\quad\qed
\enddemo

\head 3. Kashiwara modules\endhead

\subhead 3.1\endsubhead
Let $\Vb$ be any integrable $\Ub$-module and let $z$ be
a formal variable. 
Set $\Ub[z^{\pm 1}]=\Ub\otimes_\AA\AA[z^{\pm 1}].$
We denote by  $\Vb\{z\}$ be the representation of $\Ub$ on the
space $\Vb\otimes_\AA\AA[z^{\pm 1}]$ such that $x_{ir}^\pm$ acts as
$x_{ir}^\pm\otimes z^r$ and $k_{ir}^\pm$ as $k_{ir}^\pm\otimes z^r$
(in particular $e_a$ acts as $e_a\otimes z^{\delta_{ap}}$ and 
$f_a$ as $f_b\otimes z^{-\delta_{ap}}$). Define 
$$\Vb^{\otimes d}\{z_1,z_2...z_d\}=\Vb\{z_1\}\otimes_\AA\Vb\{z_2\}\otimes_\AA
\cdots
\otimes_\AA \Vb\{z_d\}.$$
The element $\Lambda\in\tilde\X$ is a weight of $\Vb$  if
the $\Lambda$-weight subspace
$$\Vb_{\Lambda}=\{v\in\Vb\,|\,
l_a\cdot v=q^{(\Lambda:\epsilon_a)}v,\,\forall a=1,...,p\}$$
is non-zero. If $v\in\Vb_{\Lambda}$ we set $\wt(v)=\Lambda.$

\subhead 3.2\endsubhead
Let $\Wb(\Lambda_a)$, $\Vb(\Lambda_a)$ be the $\AA$-form of the
fundamental $\Ub_\KK$-module and of  Kashiwara's maximal integrable
$\Ub_\KK$-module, respectively,  with highest weight $\Lambda_a$,
see \cite{K1}, \cite{K2}.
Fix  highest weight vectors $w_{\Lambda_a}\in\Wb(\Lambda_a)$,
$v_{\Lambda_a}\in\Vb(\Lambda_a)$.
By \cite{K2, Theorem 5.15(viii)} there exists a $\Ub$-automorphism $z$
of $\Vb(\Lambda_a)$ and an isomorphism of  $\Ub[z^{\pm 1}]$-modules
$$\Vb(\Lambda_a)\simto \Wb(\Lambda_a)[z]
\leqno(3.2.1)$$
such that   $v_{\Lambda_a}\mapsto w_{\Lambda_a}\otimes 1.$
In particular we have the $\Ub[z_1^{\pm 1},...,z_d^{\pm 1}]$-isomorphism
$$\Vb(\Lambda_1)^{\otimes d}\simeq
\Wb(\Lambda_1)^{\otimes d}\{z_1,z_2,...,z_d\}.$$
Both modules $\Vb(\Lambda_a),\Wb(\Lambda_a)$ have canonical bases 
$\Bb(\Lambda_a),\Bb^f(\Lambda_a),$ see \cite{K1, Proposition 8.2.2}
and \cite{K2, Theorem 5.15(ii)}.
There is a unique $\Ub$-semilinear (with respect
to $x\mapsto \overline x$) and $\QQ[z^{\pm 1}]$-linear involution
on $\Vb(\Lambda_a)$ which fixes $v_{\Lambda_a}$.
The elements in $\Bb(\Lambda_a)$ are fixed by this involution.  
Via (3.2.1), we have
$$\Bb(\Lambda_a)= \bigsqcup_{n\in\ZZ}\Bb^f(\Lambda_a)\otimes z^n.
\leqno(3.2.2)$$
Let $L(\Lambda_a)$ be the $\AA_0$-lattice spanned by $\Bb(\Lambda_a).$
\subhead 3.3\endsubhead
Let $\Vb$ be the vectorial representation of $\Ub$, i.e. $\Vb$ is the
$\Ub[z^{\pm 1}]$-module with $\AA$-basis $\{u_m\,|\,m\in\ZZ\}$ 
such that $u_{m-p\ell}=u_m\cdot z^\ell$
for all $m,\ell\in\ZZ$, and such that
$\Ub$ acts as follows ($a=1,2,...,p$)
$$
e_a(u_m)=\delta(m\equiv a+1)u_{m-1},\quad
f_a(u_m)=\delta(m\equiv a)u_{m+1},
\quad l_a(u_m)= q^{\delta(m\equiv a)}u_m.$$
Here again we write $a\equiv b$ for $a-b\in p\ZZ$. 
The $\Ub$-module $\Vb$ has a $\tilde\X^a$-graduation such that
$\wta(u_{m-pl})=\epsilon_m+l\delta$, for all $m=1,2,...,p.$
\proclaim{Lemma} 
\roster
\item In $\Vb$ we have $h_{1,1}(u_1)=o(1)(-1)^{1-p}q^{-p}u_1\cdot z.$
\item There is a unique isomorphism of $\Ub[z^{\pm 1}]$-modules 
$\Vb(\Lambda_1)\to \Vb$ such that 
$v_{\Lambda_1}\mapsto u_1.$
\item  $\Bb(\Lambda_1)
=\{u_m\,|\, m\in\ZZ\}.$
\endroster\endproclaim
\demo{Proof} For Claim 1, recall that $h_{1,1}=k_1^{-1}[x^{+}_{1,1},f_1].$
Then use a computation as in \cite{VV1, Theorem 3.3} (but note that we
use a different coproduct and different  Drinfeld generators).
 Claim (2) follows from Claim (1) and from
\cite{N, Proposition 3.1}. Claim (3) follows from (3.2.2) and standard
properties of canonical bases.
\quad\qed
\enddemo

From now on  we identify $\Vb(\Lambda_1)$ with $\Vb.$
For  $v=\otimes_iv_i\in\Vb^{\otimes d}$
and $\g=\Sum_im_i\eps_i\in\X,$ put 
$v\cdot z_\g=\otimes_i(v_i\cdot z^{m_i}).$ 
Let $\{u_\g\,|\,\g\in\X\}$ be the unique $\AA$-basis
of $\Vb^{\otimes d}$ such that
$$u_{\g-p\nu}=(\Otimes_{i=1}^du_{m_i})\cdot z_{\nu},\qquad \g,\nu\in\X.$$
Note that, if $\g\in\X_p,\nu\in\X$,
then  $\wt(u_{\g-p\nu})=\tilde\g.$
The set of  weights of $\Vb^{\otimes d}$ is exactly $\Omega.$
We have a $\tilde\X^a$-graduation on $\Vb^{\otimes d}$
such that $\wta(u_{\g-p\nu})=\tilde\g+
\sum_{i=1}^d n_i\delta$ if $\nu=\sum_{i=1}^d n_i\eps_i.$
\subhead 3.4\endsubhead 
Consider the  right representation of the $\AA$-algebra $\Hb$ on 
$\Vb^{\otimes d}$ 
such that for all $i=1,2,...d-1$,$\g\in\X_p,$
$\kappa,\nu\in\X,$
$$u_\g\diamond t_i=
\left\{\matrix
qu_\g\hfill\ &\text{if}\ \g\cdot s_i=\g\hfill\cr
u_{\g\cdot s_i}\hfill\ &\text{if}\ \g\cdot s_i<\g\hfill\cr
u_{\g\cdot s_i}+(q-q^{-1})u_\g\hfill\ 
&\text{if}\ \g\cdot s_i>\g\hfill
\endmatrix\right.,\leqno(3.4.1)$$
$$u_\nu\diamond  x_\kappa=u_{\nu\cdot\tau_\kappa}.$$
This action commutes with the $\Ub$-action, see \cite{VV2} for instance.
Set $\Tb_d=\bigoplus_\gb\rho_\gb\Hb$ where
$\gb$ runs into the set of the compositions of $d$ with at most $p$ parts.
It is a left $\Sb_d$-module such that
$t\star h=\delta(\fb',\gb)m_\gb^{-1} th\in\rho_\fb\Hb$
for all $t\in\Hb_{\fb,\fb'}$ and $h\in\rho_\gb \Hb.$
Define an involution $\iota_S$ on $\Sb_d$ by setting $\iota_S(a)=
q^{2\nu_{\fb'}}\overline {a}$, for all $a\in \Hb_{\fb,\fb'}.$

\proclaim{Lemma}
\roster
\item 
There is a 
unique isomorphism of right $\Hb$-modules
$$\phi: \Vb^{\otimes d}\simto\Tb_d,\quad
u_\mu\mapsto\rho_\eb,\qquad
\forall\, \mu\in\bar{A}_+\cap\X_p.$$
In particular 
$$\bigl(\Vb^{\otimes d}\bigr)_{\tilde\mu} \simto\rho_\eb\Hb.
\leqno(3.4.2) $$
\item There exists a unique $\AA$-algebra homomorphism $\Phi_d:\dot{\Ub}
\to \Sb_d$ such that $\phi(uv)=\Phi_d(u)\star \phi(v),$ for all 
$u\in\dot{\Ub},v\in\Vb^{\otimes d}.$ We have $\iota_S(\Phi_d(u))=
\Phi_d(\overline{u}).$
\item $u_{\mu}\diamond t_w=u_{\mu\cdot w}$
for all $\mu\in\bar A_+\cap\X$ and $w\in{}^\eb\W.$
\endroster
\endproclaim
\demo{Proof} For Claims (1) and (2) use  \cite{VV2, Lemmas 8.3 and 8.4}.
Note that the coproduct
in {\it loc.cit.} is different from our.

Claim (3) is standard, 
we prove it by induction on the lenght of $w.$
Set $w=s_i$.
Since $\mu\in \bar A_+,$ we have
$(\mu: \a_i^\vee)\geq 0,$ and the equality holds
if and only if $s_i\in\W_\eb.$ Then $\mu\cdot s_i< \mu$ and
$u_\mu\diamond t_i=u_{\mu\cdot s_i}.$
Suppose now that $\ell(w)=k>1$ and write $w=w's_i$ with
$\ell(w')=k-1.$ Then $w'\in{}^\eb\W.$
Otherwise there exists $\sigma\in \W_\eb$ such that 
$\ell(\sigma w')<\ell(w')$ and then 
$\ell(w)<\ell(\sigma w)=\ell(\sigma w')\pm 1\leq \ell(w'),$
which is impossible, because $w\in{}^\eb\W.$
So, by induction, $u_\mu\diamond t_w=
u_{\mu\cdot w'}\diamond t_i.$
Since $w'(\a_i)\in R^+,$ we have 
$(\mu\cdot w': \a_i^\vee)= (\mu:w'( \a_i^\vee) )\geq 0.$
If $(\mu\cdot w':\a_i^\vee)=0$, then
$s_{w'(\a_i)}\in\W_\eb$.
Thus $w'=s_{w'(\a_i)}w$ is an element of $\W_\eb w$  such that 
$\ell(w')<\ell(w)$ which is impossible by assumption.
\quad\qed
\enddemo

Note that the map $\Phi_d$ defined in the lemma above is a renormalization
of the map $\chi_d$ in 2.6. It is easy to check that the Claim  in section
2.6 still hold for $\Phi_{\fb,\fb'},$ the composition of $\Phi_d$ with the 
canonical projection $\Sb_d\to \Hb_{\fb,\fb'}.$
\subhead 3.5\endsubhead
Taking $d=1$ in (3.4.1) we get   
$v_{\Lambda_1}\diamond x_1= v_{\Lambda_1}\cdot z,$ because $v_{\Lambda_1}=
u_1.$ So  $x_1$ acts
as Kashiwara's operator $z$ on $\Vb.$
The module $\Wb(\Lambda_1)$ is identified with the subspace of $\Vb$
spanned by $u_1,...,u_p.$ Then, formulas (3.4.1) endow $\Wb(\Lambda_1)$
with a right $\Hb^f$-action, such that, if $p>d$,
$\Wb(\Lambda_1)^{\otimes d}
\otimes_{\Hb^f}\AA^-$ is a simple module of the quantum enveloping algebra
of $\glp$ with highest weight $\Lambda_d.$ 
Set 
$$\Vb^{[d]}=\Wb(\Lambda_1)^{\otimes d}\{(-q)^{1-d}z_1,(-q)^{3-d},\cdots,
(-q)^{d-1}z_d\}.$$
Take $\cb=(d),$ then $\Rb_\cb=\AA[\psi(x_1)].$ We can endow 
$\Vb^{[d]}
\otimes_{\Hb'}\AA^-$ with a (right)-$\Rb_\cb$-action by setting
$$(v\otimes 1)\cdot \psi(x_1)=(v\diamond x_1)\otimes 1,$$
for all $v\in\ \Vb^{[d]}.$
It is well defined by (2.1.1), (2.5.1).

\proclaim{Proposition}
Assume that $p> d$ and let $\tilde\l=\Lambda_d$.
Then $\cb=(d),$ $\rho_\db=1.$ 
Set $x=\psi((-q)^{1-d}x_1).$ 
\roster
\item
There is a unique isomorphism of $\Ub$-modules
$$\Vb(\Lambda_d)\simto\Vb^{[d]}\otimes_{\Hb'}\AA^-\quad
\text{s.t.}\quad
v_{\Lambda_d}\cdot z^{i}\mapsto (u_d\otimes \cdots u_2\otimes u_1\otimes 1)
\cdot x^i,\quad\forall i\in\ZZ.
\leqno(3.5.1)$$
\item
The set of the weights of $\Vb(\Lambda_d)$ is $\Omega_\sm.$
Fix $\tilde\mu\in\Omega_\sm.$ 
The map (3.5.1) takes $\Bb(\Lambda_d)_{\tilde\mu}$ to $\{(u_\mu\otimes 1)
\cdot x^i\,|\, i\in\ZZ\}.$ 
\endroster
\endproclaim
\demo{Proof} 
Claim (1). There exists $\a\in\AA^\times$ such that 
$(\Wb(\Lambda_1)^{\otimes d}\otimes_{\Hb^f}\AA^-)\{\a z\}\simeq
\Vb(\Lambda_d),$ because
$\Wb(\Lambda_1)^{\otimes d}\otimes_{\Hb^f}\AA^-$
and $\Wb(\Lambda_d)$ are both integral $\AA$-forms of simple $\Ub_\KK$-modules
which are simple and of
highest weight $\Lambda_d$ as modules for the quantum enveloping algebra of
$\glp.$
By (2.1.1) we have $\Rb_\cb=\AA[x^{\pm 1}].$ Consider the  $\Ub$-module  
$\Wb(\Lambda_1)^{\otimes d}\{\a z_1,q^2\a z_2,...q^{2(d-1)}\a z_d\}
\otimes_{\Hb'}\AA^-.$
The element $f_p\in\Ub$ 
acts  as
$$ \sum_{i=1}^dk_p^{i-1}\otimes\bigl(f_p\diamond x_i\a q^{2(i-1)}\bigr)
\otimes 1^{d-i}=
\sum_{i=1}^dk_p^{i-1}\otimes\bigl(f_p\diamond \a x_1\bigr)
\otimes 1^{d-i},$$
by  formula (2.5.1). A similar formula holds for $e_p.$
We need  to compute $\a$. Set $A=o(d)(-1)^{1-p}q^{-p}.$ 
Using a computation as in \cite{VV1, Theorem 3.3} we get
$$
h_{d,1}(u_d\otimes
\cdots\otimes u_1\otimes 1)=
(-q)^{d-1}Au_d\otimes\cdots\otimes u_1\diamond x_1\otimes 1.$$
Then use  \cite{N, Proposition 3.1} to get $\a=(-q)^{1-d}$. 

Claim (2). If $\mu\in(\bar A_+\setminus A_+)\cap\X_p,$
then the element $u_\mu\otimes 1\in\Vb^{[d]}
\otimes_{\Hb'}\AA^-$ is zero since $q u_\mu\otimes 1=u_\mu\diamond t_i
\otimes 1= u_\mu\otimes (t_i\cdot 1)=-q^{-1}u_\mu\otimes 1,$
for any $i\in \{1,2,...d-1\}$ such that $\mu_i=\mu_{i+1}.$ Thus 
$\Omega_\sm$ is the set of weights of $\Vb(\Lambda_d).$
Since  $f_a^2$ acts trivially on $\Vb(\Lambda_d),$ 
 Kashiwara's operator $\tilde{f_a}$ coincide with $f_a,$ for
all $a.$
Hence $\Bb(\Lambda_d)+qL(\Lambda_d)\ni
f_{a_1}\cdots f_{a_r}(v_{\Lambda_d})=u_\mu$ or $0$ with 
$\mu\in A_+\cap\X_p,$ for any $a_1,...,a_r\not=p.$ We are done by (3.2.2).
\quad\qed
\enddemo
\subhead 3.6\endsubhead
Set 
$$\Vb_\cb=\bigotimes_{i=1}^\ell\Vb(\Lambda_{c_i})\simeq
\bigotimes_{i=1}^\ell
\Wb(\Lambda_{c_i})\{z_i\},\qquad
v_\cb=\otimes_{i=1}^\ell
v_{\Lambda_{c_i}}\in\Vb_\cb.$$

It is convenient to identify the ring $\AA[z_1,...,z_\ell]$
with  the quotient of the polynomial algebra $\AA[z_1,...,z_d]$ by the ideal
generated by the relations 
$z_{\a_i}=1,$ for all $\a_i\in I_\cb.$
Let $z_{\g+\ZZ I_\cb}$ be the image of $z_\g$ in
this quotient.
Using  Proposition 3.5
we get a morphism of $\Ub$-modules
$$\bigotimes_{i=1}^\ell\Vb^{[c_i]}\to
\Vb_\cb,\quad\text {such that}\quad
\otimes_{b=p}^1\bigl(\otimes_{a=b}^1 u_a\bigr)^{\otimes\ell_b}
\mapsto v_\cb.  
\leqno(3.6.1)$$

Define a right action of the $\AA$-algebra $\Rb_\cb$ on $\bigotimes_{i=1}^\ell
\Vb^{[c_i]}
\otimes_{\Hb'_\cb}\AA^-$ by setting
$$(v\otimes 1)\cdot \psi(x_\g)=(v\diamond x_\g)\otimes 1,$$
for all $v\in\bigotimes \Vb^{[c_i]},\g\in\X.$
It is well defined by (2.1.1), (2.5.1). Via (3.4.2), we get also,
for all $\eb$,  a  right action of  $\Rb_\cb$
on $\rho_\eb\Hb\otimes_{\Hb'_\cb}
\AA^-.$ 

\proclaim{Corollary} Assume that $p>d$. 
\roster
\item
The map  (3.6.1) induces an isomorphism of  
$\Ub$-modules 
$$\Vb_\cb\simto\bigotimes_{i=1}^\ell\Vb^{[c_i]}
\otimes_{\Hb'_\cb}\AA^-,\leqno(3.6.2)$$
such that
$v_\cb\cdot z_{\g+\ZZ I_\cb}\mapsto 
(\otimes_{b=p}^1\bigl(\otimes_{a=b}^1 u_a\bigr)^{\otimes\ell_b}\otimes 1)
\cdot\psi((-q)^{-(\g:\a_\cb)}x_\g).$
\item By composing (3.6.2) with (3.4.2) we get, 
 for any $\tilde\mu,$ a unique isomorphism of $\Rb_\cb$-modules 
$$s_{\tilde\mu}: \Vb_{\cb,\tilde\mu} \simto
\rho_\eb\Hb\otimes_{\Hb'_\cb}\AA^-.$$
such that
$s_{\tilde\mu}(u\, v)=\Phi_{\eb,\fb}(u)\star s_{\tilde\nu}(v)$
if $v\in\Vb_{\cb,\tilde\nu},u\in\dot{\Ub},$ 
$uv\in\Vb_{\cb,\tilde\mu},$ and $\fb$ is
the partition associated to $\nu$ as in (1.4.1).
\endroster
\endproclaim

\subhead 3.7\endsubhead
Let $\Nb(\tilde\l)$ $\subset\Vb_\cb$ 
be the smallest $\Ub\times\Rb_\cb$-submodule containing the
element $v_\cb.$
For any composition $\fb$ of $d$, let ${}^{\fb}\W^\cb$ be the
set of elements $w\in\W^\cb$
having a maximal length in $\W_\fb w$.
The set ${}^\db\W^\cb$ contains only one element, denoted
by $\sigma_\cb,$ see \cite{L5, 13.11$(a)$}. 

For any $\mu\in \bar A_+\cap\X_p$ consider the morphism of 
$\Hb\times\Rb_\cb$-modules 
$$a_{\tilde\mu}: \rho_\eb\Hb\otimes_{\Hb'_\cb}\AA^-\to
\bigl(\bigotimes _{i=1}^\ell\Vb^{[c_i]}\otimes_{\Hb'_\cb}\AA^-\bigr)_{\tilde\mu}
$$ such that 
$(\rho_\eb h\otimes 1)\cdot \psi(x_\gamma)\mapsto (u_\mu\diamond\overline{h}
\otimes 1)\cdot \psi(x_\gamma).$

\proclaim{Proposition}
Assume that $p> d$.
\roster
\item For any weight $\tilde\mu$ we have the following
commutative square of $\Rb_\cb$-modules
$$\matrix
\Hb_{\eb,\db}
\overline{t}_{\sigma_\cb}\Rb_\cb
\otimes 1
&\hookrightarrow&
\rho_\eb\Hb\otimes_{\Hb'_\cb}\AA^-\cr
\downarrow&&\downarrow\cr
\Nb(\tilde\l)_{\tilde\mu}
&\hookrightarrow&
\Vb_{\cb,\tilde\mu}
\endmatrix$$
The vertical maps  are  $a_{\tilde\mu}$  and
its restriction.
If $\tilde \mu=\tilde \l$ then
$$a_{\tilde\l}(q^{\nu_\db}\rho_\db\overline{t}_{\sigma_\cb}
\otimes 1)=v_\cb.$$
\item  If $\tilde\mu\in\Omega_\sm$,
the left vertical map is an isomorphism. 
\endroster
\endproclaim
\demo{Proof} Claim (1).
 Take $\tilde \mu=\tilde\l$, then
$\Otimes_{a=p}^1 u_a^{\otimes d_a}=u_\l.$ 
Set $\sigma_\cb=w_\db \sigma,$ where $\sigma\in {}^\db\W.$ 
Then $t_{\sigma_\cb}=t_{w_\db} t_\sigma.$
Using (3.4.1) and Lemma 3.4(3) we get 
$$\matrix
a_{\tilde\l}(q^{\nu_\db}\rho_\db\overline{t}_{\sigma_\cb}\otimes 1)&=&
u_\l\diamond(q^{-\nu_\db}t_{\sigma_\cb})\otimes 1)\cr
&=&u_\l\diamond t_\sigma\otimes 1\cr
&=&
\otimes_{b=p}^1\bigl(\otimes_{a=b}^1 u_a\bigr)
^{\otimes\ell_b}\otimes 1.
\endmatrix
$$
Moreover 
$$\Nb(\tilde\l)_{\tilde\mu}=\bigl((\Ub\otimes
\Rb_\cb) v_{\tilde\l}\bigr)_{\tilde\mu}
{\buildrel s_{\tilde\mu}\over \lra}\,\Phi_{\eb,\db}(\dot{\Ub})\star
s_{\tilde\l}(v_\cb)\Rb_\cb.$$
If $\tilde\mu$ is small, we have $\Phi_{\eb,\db}(
\dot{\Ub})=\Hb_{\eb,\db}=\Hb\rho_\db$, by Claim 2.6 and the equality 
$\rho_\eb=1.$ Claim (2) easily follows.
\quad\qed
\enddemo
\head 4. Periodic modules\endhead

\subhead 4.1\endsubhead
Set $V':=Q\otimes_\ZZ\RR,$
 $F_{\a^\vee,n}=\{\g\in V'\,|\,(\a^\vee:\g)=np\}$
for all $\a\in R^+$, $n\in\ZZ$.
The connected components of the open set
$V'\setminus\bigcup_{\a\in R^+,n\in\ZZ}F_{\a^\vee,n}$
are called alcoves.
The connected components of the open set
$V'\setminus\bigcup_{i<d,n\in\ZZ}F_{\a_i^\vee,n}$
are called boxes.
We denote by $\Ac$ the set of alcoves.
Consider the restriction to $V'$ of the right action of $\W'\subset \W$ 
on $V.$
It induces a right $\W'$-action on $\Ac$,
denoted by $A\mapsto A\cdot w,$
which is simply transitive.
Let $A'_+$ be the unique alcove contained in the dominant Weyl chamber
having 0 in its closure.
Thus, $A'_+=A_+\cap V'.$
We define a left $\W'$-action on $\Ac$ by
$$w\cdot (A'_+\cdot v)=A'_+\cdot wv,\qquad
\forall v,w\in\W'.$$
We denote by $\leq$ the partial order on $\Ac$
defined in \cite{L1, 1.5}.
\subhead 4.2\endsubhead
We now consider the parabolic case.
Fix $\cb$  as in Convention 1.4.
The connected components of the open set
$V'\setminus\bigcup_{\a\in R^+\cap I_\cb, n\in\ZZ}F_{\a^\vee,n}$
are called $I_\cb$-alcoves.
Let $S_\cb$ be the unique $I_\cb$-alcove containing $A'_+$.
Set $\Ac_\cb\subseteq\Ac$ equal to the set of alcoves
contained in $S_\cb$.

Let $\Mb_\cb$ be the free $\AA$-module on $\Ac_\cb$.
There is a unique left action of $\Hb'$ on $\Mb_\cb$ such that
$$t_i\cdot A=\left\{
\matrix
-q^{-1} A\hfill\quad&
\text{if}\ s_i\cdot A\notin\Ac_\cb\hfill\cr
s_i\cdot A\hfill\quad&
\text{if}\ s_i\cdot A> A,\,s_i\cdot A\in\Ac_\cb\hfill\cr
s_i\cdot A+(q-q^{-1})A\hfill\quad&
\text{if}\ s_i\cdot A< A,\,s_i\cdot A\in\Ac_\cb,\hfill
\endmatrix\right.$$
for all $A\in\Ac_\cb$, $i=1,2,...d$, see \cite{L5, 9.3.(b)}.

By \cite{L5, Lemma 9.6(b)} we have $A'_+\cdot w=\overline{t}_w\cdot A'_+$,
for all $w\in\W^\cb.$ 
For any $\g\in Q$ 
there is a unique element $w_\g\in\W'_\cb$ such that
$S_\cb-p\g=S_\cb\cdot w_\g,$ because $\W'_\cb$ acts simply transitively
on the set of $I_\cb$-alcoves.
Thus the map $ A\mapsto (A-p\g)\cdot w_\g^{-1}$  
permutes the alcoves in $\Ac_\cb.$
Let  $g_\g : \Mb_\cb\to\Mb_\cb$ be corresponding $\AA$-linear map.
Then
$$g_\g(A'_+)=(-q)^{-(\g:\a_\cb)} x_\g\cdot A'_+,\quad 
g_\g(h' A'_+)=h'g_\g(A'_+),\qquad \forall h'\in\Hb',$$
see \cite{L5, 9.2 and Lemma 9.5}.
Consider the  unique right  
$\Rb'_\cb$-action on $\Mb_\cb$ such that
$$A\cdot x_{\g+\ZZ I_\cb}=(-1)^{(\g:\a_\cb)}g_\g(A),\quad
\forall A\in\Ac'_\cb,\,\forall\g\in Q$$
see \cite{L5, 9.7}. 
\subhead 4.3\endsubhead
The element
$$m_\cb=q^{-\nu_\db}\Sum_{w\in\W_\db}q^{\ell(w)}
A'_+\cdot w\sigma_\cb,$$
belongs to $\Mb_\cb$ because $\W_\db \sigma_\cb\subset\W^\cb$,
see \cite{L5, Lemma 13.9.$(b)$}.
Let $\Mb'_\cb\subset\Mb_\cb$ be the
$\Hb'\otimes\Rb'_\cb$-submodule generated by $m_\cb.$
Since, for all $w\in\W_\db$, 
$\ell(\sigma_\cb)=\ell(w\sigma_\cb)+\ell(w)$ 
we get
$$
m_\cb=
q^{-\nu_\db}\sum_{w\in\W_\db}q^{\ell(w)}\bar t_{w\sigma_\cb}A'_+
=q^{-\nu_\db}\rho_\db \bar {t}_{\sigma_\cb}\cdot A'_+.
$$
The right $\Rb_\cb$-action on $\Hb\otimes_{\Hb'_\cb}\AA^-$ defined in 3.6
restricts
to a right $\Rb'_\cb$-action on $\Hb'\otimes_{\Hb'_\cb}\AA^-.$ 

\proclaim{Lemma}
\roster
\item $\{g_\g(A'_+\cdot w)\,|\,\g\in Q, w\in \W^\cb\}
=\Ac_\cb.$
\item The linear map
$$b:\Mb_\cb\to
\Hb'\otimes_{\Hb'_\cb}\AA^-\qquad
(A'_+\cdot w)\cdot x_{\g+\ZZ I_\cb}\mapsto 
(\overline{t}_w\otimes 1)\cdot\psi(q^{-(\g:\a_\cb)}x_\g),$$
where $w\in\W^\cb, \g\in Q$,
is a morphism of $\Hb'\times\Rb'_\cb$-modules such that
$b(m_\cb)=q^{-\nu_\db}\rho_\db
\bar{ t}_{\sigma_\cb}\otimes 1.$ 
\item The map $b$ yields a  commutative square of
$\Hb'\otimes\Rb'_\cb$-modules
$$\matrix
\Hb'\rho_\db \bar t_{\sigma_\cb}\Rb'_\cb\otimes 1
&\hookrightarrow&
\Hb'\otimes_{\Hb'_\cb}\AA^-\cr
\uparrow&&\uparrow\cr
\Mb'_\cb
&\hookrightarrow&
\Mb_\cb
\endmatrix$$
with invertible vertical maps.  
\endroster
\endproclaim
\demo{Proof}
Claim (1) follows from \cite{L3, 2.12$(f)$, 4.9$(c)$} and the identity 
$w_d\W^\cb w_\cb=\W^\cb.$ 
As for  Claims (2) and (3), use 
\cite{L5, 9.7$(b)$ and Proposition 8.5, Lemma 10.2}.
\quad\qed
\enddemo
\head 5. Comparing canonical bases\endhead

\subhead 5.1\endsubhead
Let  $\Mb_{\cb,\leq}$ be the set of formal $\AA$-linear combinations
$m=\sum_{A\in\Ac_\cb}^\infty m_A\,A$ such that the set 
$\{A\,|\,m_A\neq 0\}$ is bounded above under $\leq,$ i.e. there 
exists $B\in\Ac_\cb$ such that $A< B$ for all $A$ such that
$m_A\neq 0.$
The $\Hb'\times\Rb'_\cb$-action on $\Mb_\cb$
extends in an obvious way to $\Mb_{\cb,\leq}$.

There is a unique continuous (in the sense of \cite{L3, 4.13})
involution
$\iota_M:\Mb_{\cb,\leq}\to\Mb_{\cb,\leq}$ such that,
for all $m\in\Mb_{\cb,\leq}, h\in\Hb',x_{\g+\ZZ I_\cb}\in\Rb'_\cb,$
$$\iota_M(m_\cb)=m_\cb,\quad
\iota_M(hm)=\overline{h}\,\iota_M(m),\quad 
\iota_M(m\cdot x_{\g+\ZZ I_\cb})=\iota_M(m)\cdot x_{\g+\ZZ I_\cb}
,\leqno(5.1.1)$$
see \cite{L5, Lemmas 9.16.$(a)$ and  13.12}.

For all $A\in\Ac_\cb$ there is a unique element
$A_\leq=\Sum^\infty_{B\preceq A}\tilde\pi_{BA}B\in\Mb_{\cb,\leq}$
such that 
$$\iota_M( A_\leq)=A_\leq,\qquad \tilde\pi_{AA}=1\and
\tilde\pi_{BA}\in q^{-1}\ZZ[q^{-1}]\text{ for all }B< A,$$
see \cite{L5,9.17(a),(b)}.
We have
$m_\cb=(A'_+\cdot w_\db\sigma_\cb)_\leq$ by
\cite{L5, Theorem 13.13$(a)$}.
Set $\Bc_M=\{A_\leq\,|\, A\in \Ac_\cb\}.$

\proclaim{Conjecture \cite{L5, 12.7}}
$\Bc_M\subset \Mb_\cb.$ 
\endproclaim

For a future use let us recall the following standard fact.
\proclaim{Lemma} Let $A^\sharp=A+\sum_{A'\not=A} m_{A'}A'
\in\Mb_\cb,$ with  $m_{A'}\in q^{-1}\QQ[q^{-1}].$
If  $\iota_M(A^\sharp)=A^\sharp,$ then $A^\sharp\in\Bc_M.$
\endproclaim
\demo{Proof} For all $B\in\Ac_\cb$ we have
$B=B_\leq+\Sum_{B'< B}^\infty c_{B'}B'_\leq$
in $\Mb_{\cb,\leq},$ with $c_{B'}\in q^{-1}\QQ[q^{-1}].$
Thus $A^\sharp=\sum_B^\infty n_B B_\leq,$ with $n_B\in q^{-1}\QQ[q^{-1}]$ for
all $B\neq A,$ and $n_A\in 1+q^{-1}\QQ[q^{-1}].$ 
Since $\iota_M(A^\sharp)=A^\sharp,$ we get 
$\overline{n}_B=n_B$ for all $B.$ Thus $A^\sharp=A_\leq.$
\quad\qed
\enddemo

\subhead 5.2\endsubhead
Following \cite{K2, Theorem 8.5, Proposition 8.6},
the $\AA$-module $\Nb(\tilde\l)$ is endowed with
a canonical basis. More precisely, define
$$\Nb_{\cb,\leq}=\Vb_\cb
\otimes_{\QQ[z_{\a_i+\ZZ I_\cb}|\a_i\in I]}
\QQ[[z_{\a_i+\ZZ I_\cb}|\a_i\in I]],$$ 
$$L_{\cb,\leq}=\bigl(\bigotimes_{i=1}^\ell 
L(\Lambda_{c_i})\bigr)
\otimes_{\QQ[z_{\a_i+\ZZ I_\cb}|\a_i\in I]}
\QQ[[z_{\a_i+\ZZ I_\cb}|\a_i\in I]],$$ 
There is a unique  involution
$\iota_N:\Nb_{\cb,\leq}\to\Nb_{\cb,\leq}$
 such that,
for all $v\in\Nb_{\cb,\leq},
u\in\Ub, f\in\QQ[[z_{\a_i+\ZZ I_\cb}|\a_i\in I]][z_\g|\g\in\X],$
$$\iota_N(v_\cb)=v_\cb,\quad
\iota_N(uv)=\overline{u}\,\iota_N(v),
\quad \iota_N(v\cdot f)= \iota_N(v)\cdot f
.\leqno(5.2.1)$$
For any $t\in\bigotimes_i\Bb(\Lambda_{c_i})$ there is a unique element
$F(t)\in \Nb_{\cb,\leq}$ such that 
$$\iota_N(F(t))=F(t),\qquad F(t)= t\quad
mod \,  qL_{\cb,\leq}.$$
The $\AA$-module $\Vb(\Lambda_{c_i})$
has a $\tilde\X^a$-graduation induced by  the $\tilde\X^a$-graduation 
of $\Vb^{\otimes d},$ using Proposition 3.5. Let $\wta(v)\in \tilde\X^a$
be the degree of the element $v\in \Vb(\Lambda_{c_i}).$
There is a unique partial order, $\leq_\cb$, 
on $\bigotimes_i\Bb(\Lambda_{c_i})$
such that $t'=t_1'\otimes\cdots\otimes t_\ell'<_\cb 
t=t_1\otimes\cdots t_\ell$ if and only if
$$\sum_{i=1}^\ell\wta(t'_i)=\sum_{i=1}^\ell\wta(t_i),\quad
\sum_{i=1}^m(\wta(t'_i)-\wta(t_i))\in 
\tilde{Q}^a_+\setminus\{0\},\quad \forall m=1,...,\ell-1,$$
where  $\tilde{Q}^a_+=\sum_{a=1}^p\NN\b_a$.
Then
$$F(t)-t\in\Sum_{t'<_\cb t}q\QQ[q]t'.\leqno(5.2.2)$$
The $\Ub\times\Rb_\cb$-submodule $\Nb(\tilde\l)$ of 
$\Nb_{\cb,\leq}$ is stable by $\iota_N.$ Furthermore 
$F(t)\in\Nb(\tilde\l).$ 
The family $\Bc_N=\{F(t)\,|\, t\in\bigotimes_i \Bb(\Lambda_{c_i})\}$ is
the canonical basis of $\Nb(\tilde\l).$
\subhead 5.3\endsubhead
The basis  $\bigotimes_{i=1}^\ell \Bb(\Lambda_{c_i})$ of $\Vb_\cb$
is identified with its image in $\bigoplus_{i=1}^\ell\Vb^{[c_i]}
\otimes_{\Hb'_\cb}\AA^-$ by (3.6.2).

\proclaim{Claim 1} $\bigotimes_{i=1}^\ell \Bb(\Lambda_{c_i})=
\{(u_{\mu\cdot w}\otimes1)\cdot z_{\g+\ZZ I_\cb}\,|\,
\mu\in A_+\cap\X_p,
w\in \W^\cb,\g\in\X\}.$
\endproclaim
\demo{Proof} 
Given $\mu\in A_+\cap\X_p,$  $w\in\W^f,$ write $\mu\cdot w=
\sum_{j=1}^da_j\eps_j.$ Then $w\in\W^\cb$ if and only if $p>a_k-a_l>0$
for all $c_1+\cdots+c_i\geq k>l\geq c_1+\cdots c_{i-1}+1.$
The  claim follows from  Proposition 3.5(2).
\quad\qed
\enddemo
Assume that $p>d,$ then $\Omega_\sm\not=\emptyset.$ 
Take $\tilde\mu\in\Omega_\sm.$
The  map  $d_{\tilde\mu}=a_{\tilde\mu}
\circ b$ gives an injective morphism of $\Hb'\times\Rb'_\cb$-modules
$$d_{\tilde\mu}:\Mb_\cb\to\bigotimes_{i=\ell}^1\Vb^{[c_i]}
\otimes _{\Hb'_{\cb}}\AA^-,$$
$$(A'_+\cdot w)\cdot x_{\g+\ZZ I_\cb}\mapsto 
(-1)^{(\g:\a_\cb)}(u_{\mu}\diamond t_w \otimes 1)\cdot z_{\g+\ZZ I_\cb},$$
for all $w\in\W^\cb,\g\in Q.$
Moreover $d_{\tilde\mu}$ maps $\Mb_\cb'$  into
$\Nb(\tilde\l)_{\tilde\mu}$ by Lemma 4.3(3) and Proposition 3.7(1).

\proclaim{Claim 2} $d_{\tilde\mu}(\Ac_\cb)\subset 
\bigotimes_{i=1}^\ell \Bb(\Lambda_{c_i}).$ 
\endproclaim
\demo{Proof} By Lemma 4.3(1) any alcove in $\Ac_\cb$ is of the
type $A=g_\g(A'_+\cdot w)$ for some $\g\in Q,$ $w\in\W^\cb$.
Thus $d_{\tilde\mu}(A)=(u_\mu\diamond t_w\otimes 1)\cdot z_{\g+\ZZ I_\cb}$
and Claim 2 follows from Claim 1 and Lemma 3.4(3) 
(remember that $\mu\in A_+\cap\X_p$).
\quad\qed
\enddemo
\subhead 5.4\endsubhead
Set $\X'=(\bar{A}_+\cap\X_p)\cdot \W'.$ 

\proclaim{Lemma} There exists a family $\Bc'_N\subset\Bc_N$ such that 
the $\AA$-span of $\Bc'_N$ is equal to the intersection of $\Nb(\tilde\l)$ with
the $\AA$-span of $\{u_\g\otimes 1
\,|\, \g\in\X'\}.$
\endproclaim
\demo{Proof} First of all, $\wta(u_{\g_1})=\wta(u_{\g_2})$ iff
$\g_1\in\g_2\cdot \W',$ 
for any $\g_1,\g_2\in \X.$ Namely,
write $\g_i=\nu_i-p\mu_i$ with $\nu_i\in\X_p,\mu_i
\in\X.$ Then $\wta(u_{\g_1})=\wta(u_{\g_2})$ iff 
$\mu_1-\mu_2\in Q$ and $\tilde\nu_1=\tilde\nu_2$ i.e. $\nu_1=\nu_2\cdot w$
for some $w\in\W^f,$ hence iff $\g_1=\g_2\cdot w\tau_{\mu_1-\mu_2\cdot w}$
and $\mu_1-\mu_2\cdot w\in Q.$
 
Given $\nu,\g\in \X$ such that $u_\nu\otimes 1,u_\g\otimes 1\in
q^\ZZ\bigotimes_{i=1}^\ell\Bb(\Lambda_{c_i}),$ we write $u_\nu\otimes 1
<_\cb u_\g\otimes 1$ if the corresponding elements in 
$\bigotimes_{i=1}^\ell\Bb(\Lambda_{c_i})$ satisfy the same relation.
Then
$u_\nu\otimes 1<_\cb u_{\g}\otimes 1$  implies $\wta(u_\nu)=\wta(u_\g),$
hence 
$\nu\in\g\cdot \W'.$ 
For any $\g\in\X,$ set $F(u_\g\otimes 1)=F(q^nu_\g\otimes 1)$ 
if $q^nu_\g\otimes 1
\in \bigotimes_{i=1}^\ell\Bb(\Lambda_{c_i}).$
We have $F(u_\g\otimes1)\in \bigoplus_{w\in\W'} \AA u_{\g\cdot w}\otimes 1$
by (5.2.2). Thus, for all $\g$,
$$\bigoplus_{w\in \W'}\AA u_{\g\cdot w}\otimes 1=
\bigoplus_{w\in \W'}\AA F(u_{\g\cdot w}\otimes 1).$$

Since $\Nb(\tilde\l)$ is spanned by the elements $F(t)$ with $t
\in \bigotimes_{i=1}^\ell\Bb(\Lambda_{c_i})$ we are done because 
$\X'\cdot \W'\subseteq \X'.$
\quad\qed
\enddemo
Fix $p>d$ and $\tilde\mu\in\Omega_\sm.$ 
Set $\Bc'_{N,\tilde\mu}=\Bc'_N\cap \Nb(\tilde\l)_{\tilde\mu},$ and let 
$\Nb'(\tilde\l)_{\tilde\mu}$ be the $\AA$-span of  $\Bc'_{N,\tilde\mu}.$
The $\AA$-span of 
$\{u_\mu\,|\,\mu\in\X'\}$ in $\Vb^{\otimes d}$ is
identified to 
$\bigoplus_\gb\rho_\gb\Hb',$ see (3.4.1) and Lemma 3.4(1).
We have  $a_{\tilde\mu}^{-1}(\Nb'(\tilde\l)
_{\tilde\mu})=\Hb'\rho_\db\overline{t}_{\sigma_\cb}\Rb'_\cb\otimes 1$
since 
$a_{\tilde\mu}^{-1}(\Nb'(\tilde\l)_{\tilde\mu})=(\rho_\eb\Hb'\otimes _{\Hb'_\cb}
\AA^-)\cap (\Hb\rho_\db\overline{t}_{\sigma_\cb}\Rb_\cb\otimes 1).$
Hence $\Nb'(\tilde\l)_{\tilde\mu}$ is a right $\Rb_\cb'$-module by the 
lemma above and (3.4.1), and a left $\Hb'$-module. 
We have $d_{\tilde\mu}(\Mb'_\cb)=a_{\tilde\mu}(\Hb'\rho_\db
\overline{t}_{\sigma_\cb}\Rb'_\cb\otimes 1)=\Nb'(\tilde\l)_{\tilde\mu}.$
Let $c_{\tilde\mu}$ be the inverse map
$c_{\tilde\mu}:
 \Nb'(\tilde\l)_{\tilde\mu}\to \Mb'_\cb.$
It is an isomorphism of $\Hb'\otimes \Rb_\cb'$-modules.

\subhead 5.5\endsubhead
We can now prove Conjecture  5.1.
\proclaim{Theorem} $\Bc_M$  is a $\AA$-basis of $\Mb'_\cb$ for all $p>d.$
\endproclaim
\demo{Proof} Let $\tilde\mu\in\Omega_\sm.$ Set $h=q^{-\nu_\db}\rho_\db
\overline{t}_{\sigma_\cb}.$ We consider $h$ as an element of
$\rho_\eb\Hb\subset \Tb_d.$ 
Define  
$v_{\cb,\tilde\mu}=a_{\tilde\mu}(h\otimes 1).$ Then
$$v_{\cb,\tilde\mu}=
u_\mu\diamond \overline{h}\otimes 1=
u_\mu\diamond q^{-\nu_\db}\rho_\db t_{\sigma_\cb}\otimes 1=
d_{\tilde\mu}(m_\cb) \in\Nb'(\tilde\l)_{\tilde\mu}.$$
The element $v_\cb=
a_{\tilde\l}(q^{\nu_\db}\rho_\db\overline{t}_{\sigma_\cb}\otimes 1)=
u_\l\diamond q^{-\nu_\db}t_{\sigma_\cb}\otimes 1\in
\Nb'(\tilde\l)_{\tilde\l}$ is fixed by $\iota_N.$ Hence
$\iota_N(v_{\cb,\tilde\mu})=
v_{\cb,\tilde\mu},$  because $v_{\cb,\tilde\mu}=\rho_\db\star
v_\cb$ and $\iota_S(\rho_\db)=\rho_\db.$

For any  $v\in \Nb'(\tilde\l)_{\tilde\mu}$ there exists $u
\in \dot{\Ub}$ and $r\in \Rb'_\cb$  such that $v=uv_{\cb,\tilde\mu}
\cdot  r,$ since $\Phi_{\eb,\db}(\dot{\Ub})\star (q^{-\nu_\db}\rho_\db
\overline{t}_{\sigma_\cb}\otimes 1)\subseteq \Phi_{\eb,\eb}(\dot{\Ub})\star
(h\otimes 1),$ by Claim 2.6.
We want to prove that $c_{\tilde\mu}$ is compatible with the involutions 
$\iota_N,\iota_M.$ By (5.1.1), (5.2.1) and the $\Rb'_\cb$-linearity of
$c_{\tilde\mu},$ we can forget $r.$ We have 
$c_{\tilde\mu}(uv_{\cb, \tilde\mu})= 
\overline{\Phi_{\eb,\eb}(u)}\cdot m_\cb.$  
Namely, set $t=\Phi_{\eb,\eb}(u).$ Then
$a_{\tilde\mu}(uv_{\cb,\tilde\mu})=
a_{\tilde\mu}(t\star (u_\mu\diamond\overline{h}\otimes 1))=
a_{\tilde\mu}(u_\mu\diamond
t\overline{h}\otimes 1)=
\overline{t}h=
b(\overline{t}\cdot m_\cb).$

Using  5.1.1, 5.2.1, and Lemma 3.4(2) we get 
$$\matrix
c_{\tilde\mu}(\iota_N(uv_{\cb,\tilde\mu}))&=&
c_{\tilde\mu}(\overline{u}v_{\cb,\tilde\mu})\cr
&=&\overline{\Phi_{\eb,\eb}(\overline{u})}\cdot m_\cb\cr
&=&
\overline{\iota_S(\Phi_{\eb,\eb}(u))}\cdot m_\cb\cr
&=&\Phi_{\eb,\eb}(u)\cdot m_\cb\cr
&=&\iota_M(\overline{\Phi_{\eb,\eb}(u)}\cdot m_\cb)\cr
&=&\iota_M(c_{\tilde\mu}(uv_{\cb,\tilde\mu})).
\endmatrix$$

Hence any element in  $c_{\tilde\mu}(\Bc'_{N,\tilde\mu})$  is fixed by 
$\iota_M.$
Claim 2 in 5.3 gives
$$c_{\tilde\mu}(\Nb'(\tilde\l)_{\tilde\mu}\cap L_{\cb,\leq})\subseteq
\sum_{A\in\Ac_\cb}\QQ[q^{-1}]A.$$
Thus,  by Lemma  5.1, 
$c_{\tilde\mu}(\Bc'_{N,\tilde\mu})\subseteq \Bc_M.$
We are done since $c_{\tilde\mu}$ is invertible .
\quad\qed
\enddemo

\noindent{\bf Remarks.} 
\roster
\item As a corollary of Theorem 5.5 and \cite {L5, 15.13},
the set $\Bb^\pm _{\Bc_e}$ in \cite{L5, 5.11} is a signed basis of the K-theory 
of the Springer fiber. 
\item 
Even in the $A_2$ case the two orders are not comparable via $d_{\tilde\mu}.$
Fix $d=3,p>3,\cb=(1^3),$
and $\mu=3\eps_1+2\eps_2+\eps_3.$ Then 
$d_{\tilde\mu}(A'_+)=u_\mu, d_{\tilde\mu}
(A'_+\cdot w)=u_{\mu\cdot w},$ for any $w\in \W^f.$

- If $w=s_2,$  we have that  $A'_+\cdot w< A'_+,$
but $u_\mu\not<_\cb u_{\mu\cdot w}$.

- If $w=s_{\theta^\vee}\tau_{-\a_1-2\a_2},$
we have $u_\mu<_\cb u_{\mu\cdot w}$ but $A'_+\not< A'_+\cdot w,$
since $A'_+\cdot w$ is not in the Weyl chamber.
\endroster


\vskip1cm
\Refs
\widestnumber\key{ABC}

\ref\key{B}\by Beck, J.\paper Braid group action and quantum
affine algebras\jour Comm. Math. Phys.\vol 165\yr 1994\pages 555-568
\endref

\ref\key{K1}\by Kashiwara, M.
\paper Crystal bases of the modified quantized enveloping algebra 
\jour Duke Math. J.\vol 73\yr 1994\pages 383-413\endref

\ref\key{K2}\by Kashiwara, M.
\paper On level zero representations of quantized affine algebras.
\jour Duke Math. J.\vol 112\yr 2002\pages 117-175\endref

\ref\key{L1}\by Lusztig, G.\paper Hecke algebras and Jantzen's generic 
decomposition patterns\jour Adv. in Math.\vol 37\yr 1980\pages
121-164\endref

\ref\key{L2}\by Lusztig, G.\book Introduction to quantum groups 
\publ Birkh\"auser\publaddr Boston-Basel-Berlin \yr 1994\endref

\ref\key{L3}\by Lusztig, G.\paper Periodic W-graphs
\jour Represent. Theory\vol 1\yr 1997\pages 207-279\endref

\ref\key{L4}\by Lusztig, G.\paper Bases in equivariant K-theory
\jour Represent. Theory\vol 2\yr 1998\pages 298-369\endref

\ref\key{L5}\by Lusztig, G.\paper Bases in equivariant K-theory, II
\jour Represent. Theory\vol 3\yr 1999\pages 281-353\endref

\ref\key{L6}\by Lusztig, G.\paper Aperiodicity in quantum affine 
$\gln$ \jour Asian J. Math.\vol 3\yr 1999\pages 147-177\endref

\ref\key{N}\by Nakajima, H.\paper Extremal weight modules of quantum
affine algebras\jour math.QA/0204183 \endref

\ref\key{SV}\by Schiffmann, O., Vasserot, E.\paper Geometric
construction of the global base of quantum modified algebra of
$\widehat{\gln}$
\jour Thansform. Groups\vol 5\yr 2000\pages 351-360\endref

\ref\key{VV1}\by Varagnolo, M., Vasserot, E.
\paper Schur duality in the toroidal setting
\jour Commun. Math. Phys.\vol 182\yr 1996\pages 469-484\endref

\ref\key{VV2}\by Varagnolo, M., Vasserot, E.
\paper On the decomposition matrices of the quantized Schur algebra
\jour Duke Math. J.\vol 100\yr 1999\pages 267-297\endref

\ref\key{VV3}\by Varagnolo, M., Vasserot, E.
\paper Canonical bases and quiver varieties
\jour Represent. Theory\vol 7\yr 2003\pages 227-258\endref

\endRefs
\enddocument